\documentclass{amsart}
\usepackage{amsmath}
\usepackage{amssymb}
\usepackage{array}
\usepackage{amscd}
\usepackage{tikz-cd}

\numberwithin{equation}{section}
\newtheorem{Th}[subsection]{Theorem}
\newtheorem*{Th*}{Theorem}
\newtheorem{Lemma}[subsection]{Lemma}

\theoremstyle{definition}
\newtheorem{definition}[subsection]{Definition}
\newtheorem*{definition*}{Definition}
\newtheorem{Remark}[subsection]{Remark}

\newcommand{\comm}[1]{}

\usepackage{color}
\usepackage[normalem]{ulem}
\definecolor{DarkGreen}{rgb}{0,0.5,0.1} 

\newcommand\soutD{\bgroup\markoverwith
{\textcolor{DarkGreen}{\rule[.5ex]{2pt}{1pt}}}\ULon}

\makeatletter
\newcommand*{\rom}[1]{\expandafter\@slowromancap\romannumeral #1@}
\makeatother

\begin{document}
 
\title[Severi--Brauer varieties of prime degree over $\mathbb{Q}$]{Automorphism groups of Severi--Brauer varieties of odd prime degree over $\mathbb{Q}$}
\author[Raisa V.~Serova]{Raisa Vadimovna Serova}
\address{{\sloppy
\parbox{0.9\textwidth}{\emph{Raisa V.~Serova}
\newline
School of Mathematics, The University of Edinburgh, Edinburgh EH9 3JZ, UK
}\bigskip}}
\email{serova.r@yandex.ru}

\author[Anastasia V.~Vikulova]{Anastasia Vadimovna Vikulova}
\address{{\sloppy
\parbox{0.9\textwidth}{\emph{Anastasia V.~Vikulova}
\newline
Department of Mathematics and Computer Science, University of Basel, Spiegelgasse
1, 4051 Basel, Switzerland
}\bigskip}}
\email{vikulovaav@gmail.com}

\date{}
\maketitle

\begin{abstract}
We prove that the only non-trivial finite subgroup of automorphism group of non-trivial Severi--Brauer variety of odd prime degree~$p$ over the field of rational numbers is~$\mathbb{Z}/p\mathbb{Z}.$
\end{abstract}

	\section{Introduction}

Nowadays \textit{equivariant} birational geometry attracts considerable interest, since it is sometimes important to keep track of the action of a fixed group on a variety when studying its birational geometry. This motivates the study of finite subgroups of automorphism groups of varieties. In this paper we mainly deal with forms of projective spaces, known as Severi--Brauer varieties. Let us recall their definition.

\begin{definition}
An $n$-dimensional variety $V$ over a field $\mathbf{F}$ is called a \textit{Severi--Brauer variety of degree $n+1$} if 
$$
V \times_{\mathrm{Spec}(\mathbf{F})}\mathrm{Spec}(\overline{\mathbf{F}}) \simeq \mathbb{P}^n_{\overline{\mathbf{F}}},
$$

\noindent where $\overline{\mathbf{F}}$ is an algebraic closure of $\mathbf{F}.$ Such a variety $V$ is called \textit{non-trivial} if it is not isomorphic to $\mathbb{P}^n_{\mathbf{F}}.$ 
\end{definition}

The following result about finite automorphism group of Severi--Brauer surface over the field~$\mathbb{Q}$ of rational numbers was proven in~\cite{VikulovaSB}.

\begin{Th}[{\cite[Theorem 1.3]{VikulovaSB}}]\label{Z/3Z}
Let $V$ be a non-trivial Severi--Brauer surface  over the field $\mathbb{Q}$ and let $G$ be a nontrivial finite subgroup in $\mathrm{Aut}(V)$. Then~\mbox{$G \simeq\mathbb{Z}/3\mathbb{Z}.$}
\end{Th}

The goal of this paper is to extend Theorem~\ref{Z/3Z} to Severi--Brauer varieties of arbitrary odd prime degree.

\begin{Th}\label{Z/pZ}
Let $V$ be a non-trivial Severi--Brauer variety of odd prime degree~$p$  over the field $\mathbb{Q}$ and let $G$ be a nontrivial finite subgroup in $\mathrm{Aut}(V)$. Then~\mbox{$G \simeq \mathbb{Z}/p\mathbb{Z}.$}
\end{Th}

\begin{Remark}
Note that Theorem~\ref{Z/pZ} fails when $p=2.$ Indeed, consider a non-trivial Severi--Brauer curve~$C$ over~$\mathbb{Q}$ defined by the equation~$x^2+y^2+z^2=0$ in~$\mathbb{P}^2.$ Then we have that symmetric group $\mathfrak{S}_3$ lies in $\mathrm{Aut}(C).$     
\end{Remark}

The plan of the paper is as follows. 
In Section~\ref{sec:preliminaries} we prove some supplementary lemmas.
In Section~\ref{section:SBV}  we collect some basic facts about Severi--Brauer varieties and study finite subgroups of their automorphism groups.
In Section~\ref{section:proofoftheorem} we prove Theorem~\ref{Z/pZ}.

\textbf{Acknowledgment.} The authors are deeply grateful to Ivan Cheltsov and Constantin Shramov for inspiration, constant support and endless stream of wild ideas. 

The first author is supported by the Simons Collaboration grant \textit{Moduli of varieties}. The second author is supported by the project SERI REF-1131-552105.

\textbf{AI disclosure.} No AI bots were harmed and used in obtaining the mathematical results of this paper.

\section{Preliminaries}\label{sec:preliminaries}

In this section we collect some auxiliary facts.

\begin{Lemma}\label{lemma:representationz/pz}
Let $p$ be an odd prime number. Then $(\mathbb{Z}/p\mathbb{Z})^2 \not\subset \mathrm{GL}_p(\mathbf{\mathbb{Q}}).$

\end{Lemma}

\begin{proof}  Assume that there is a linear action of $(\mathbb{Z}/p\mathbb{Z})^2$ on the vector space $\mathbb{Q}^p.$ As this group is abelian we get that the matrices, which represent the elements of the group, are simultaneously diagonalizable  over~$\overline{\mathbb{Q}}$. Thus, the elements of the group are conjugate to
\begin{equation}\label{eq:matricesomegaabcdwitout323}
\begin{pmatrix}
\zeta_p^{a_0} & 0 & \cdots  & 0\\
0 & \zeta_p^{a_1} & \cdots & 0 \\
\vdots & \vdots & \ddots & 0\\
0 & 0 & \cdots & \zeta_p^{a_{p - 1}} 
\end{pmatrix},
\end{equation}

\noindent where $\zeta_p$ is a non-trivial $p$-th root of unity and $a_i \in \{0,1, \cdots, p - 1\}.$ Note that the determinant and the trace of these matrices belong to $\mathbb{Q}.$ Therefore, such matrices have to satisfy the following:
\begin{gather}
a_0 + \cdots + a_{p - 1} \equiv 0 \mod p; \label{eq:123} \\
\zeta_p^{a_{0}}+\cdots+\zeta_p^{a_{p - 1}} \in \mathbb{Q}. \label{eq:223}
\end{gather}
\noindent Let $(a_0, a_1, \ldots, a_{p-1})$ be a tuple associated with the matrix~\eqref{eq:matricesomegaabcdwitout323}.   There are~\mbox{$p! + 1$} tuples  satisfying~\eqref{eq:123} and~\eqref{eq:223}, namely, $(0, 0, \cdots, 0)$ and all permutations of the tuple~$(0, 1, 2, \cdots, p - 1)$. Let $\alpha$ and $\beta$ be two generators of~$(\mathbb{Z}/p\mathbb{Z})^2.$
Without loss of generality, we can assume that $\alpha$ and $\beta$ are associated with the tuples~\mbox{$(0, 1, \cdots, p - 1)$} and $(a_0, a_1, \cdots, a_{p-1})$, respectively, such that  one matrix is not the power of the other. Multiplying these two matrices we get the matrix~\eqref{eq:matricesomegaabcdwitout323} with the tuple
$$
(a_{0}, 1 + a_1, \cdots, p - 1 +a_{p - 1}).
$$
Assume that $a_0 \neq 0$, then $a_i = 0$ for some $i$ and the result of multiplication is the following:
$$(a_{0}, 1 + a_1, \cdots, i ,\cdots,  p - 1 +a_{p - 1}).$$
Keep multiplying the result by matrix with tuple $(a_0, a_1, \cdots, a_{p-1})$ we get
$$(i, a_1^{\prime}, \cdots, i, \cdots, a_1^{\prime}).$$
This matrix does not satisfy ~\eqref{eq:223}.
If $a_0 = 0$, acting similarly we get a matrix with tuple, where there are two zeros, which also does not satisfy ~\eqref{eq:223}.
 Hence, there is either trivial group or group of order $p$, generated by the matrix with tuple 
 $$
 (a_0, a_1, \cdots, a_{p-1}).
 $$

\end{proof}

\begin{Lemma}\label{lemma:fielddegree3}
 Let us consider the field extension $\mathbb{Q} \subset \mathbf{L}$ of odd prime degree $p$ generated by an element $\bar{x} \in \mathbf{L}$ such that $\bar{x}^p-\alpha=0$ for $\alpha \in \mathbb{Q}^*.$ Let us consider an element~\mbox{$\bar{y} \in \mathbf{L}$} such that $\bar{y}^p-\beta=0$ for $\beta \in \mathbb{Q}^*.$  Then  $\bar{y}=c\bar{x}^i$  for $c \in \mathbb{Q}^*$ and some integer~\mbox{$0 \leqslant i \leqslant p-1.$}
\end{Lemma}

\begin{proof}
First of all, assume that $\beta \in \left(\mathbb{Q}^*\right)^p.$ Then  $\bar{y}$ is an element in  $\mathbb{Q}.$ Therefore, we can assume that $\beta \notin  \left(\mathbb{Q}^*\right)^p,$ thus the element $\bar{y}$ generates a field extension of degree $p.$ This means that $1$, $\bar{y}, \cdots, \bar{y}^{p - 1}$ is a basis of $\mathbf{L}.$

We have
$$
\bar{y}=u_{0}+u_{1}\bar{x}+ \cdots +u_{p - 1}\bar{x}^{p - 1}
$$
\noindent for $u_i\in \mathbb{Q}.$ Let us consider the trace of the multiplication by the element $\bar{y}.$  Note that the trace of all elements $\bar{x}^{i}$ is zero, which can be easily computed in the basis~$1$,~$\bar{x}, \cdots ,\bar{x}^{p - 1}.$ Therefore, on the one hand, the trace of multiplication by $\bar{y}$ is equal to $pu_0.$ On the other hand, if we consider the multiplication by the element~$\bar{y}$ in the basis $1$, $\bar{y}, \cdots, \bar{y}^{p - 1},$ we get that its trace is equal to zero. So we obtain  $$\bar{y}=u_{1}\bar{x}+ \cdots +u_{p - 1}\bar{x}^{p - 1}.$$ 
Let $i$ be the smallest index such that $u_i \neq 0$, then we can write $\bar{y} = \bar{x}^i\bar{z}$, where 
$$
\bar{z} = u_i + u_{i + 1}\bar{x} + \cdots + u_{p - 1}\bar{x}^{p - 1 - i}.
$$ 
If $u_j = 0$ for all $j \neq i$, then $\bar{z} = u_i$ and we are done. Assume there exist a coefficient~\mbox{$u_j \neq 0$} with $j>i$, then we can apply the same approach as to $\bar{y}$ to $\bar{z}$. Indeed, consider the trace of the multiplication by the element $\bar{z}$. On the one hand it is zero, on the other hand it is $pu_i$ and $u_i = 0$, which is a contradiction. So we get $\bar{y} = u_i\bar{x}^i,$ as desired.

\end{proof}

\section{Severi--Brauer varieties}\label{section:SBV}
In this section we study finite subgroups of the automorphism groups of  Severi--Brauer varieties. Let us mention some properties of Severi--Brauer varieties (for more details see, for instance,~\cite{SimpleAlgebras} and~\cite{Kollar}). If $V$ is a Severi--Brauer variety over a field~$\mathbf{F},$ then $V$ is non-trivial  if and only if~\mbox{$V(\mathbf{F})= \varnothing$} (see~\cite[Theorem 5.1.3]{SimpleAlgebras}). There is a bijection between Severi--Brauer varieties of dimension~$n$ over a field~$\mathbf{F}$ and central simple algebras of dimension~$(n+1)^2$ over  $\mathbf{F}$  (see e.g.~\cite[\S 5.2]{SimpleAlgebras}).

\begin{Th}[{\cite[Theorems 2.1.3 and 5.2.1]{SimpleAlgebras}}]\label{CSAcorps}
Let $V$ be a non-trivial Severi--Brauer variety of prime degree $p.$  Then the central simple algebra $A$ which is associated to~$V$ is a division algebra.  
\end{Th}

The following theorem describes automorphism groups of Severi--Brauer varieties.

\begin{Th}[{see, for example,~\cite[p.~266]{Chatelet} or \cite[Lemma 4.1]{ShramovVologodsky}}]\label{csaprop}
If a central simple algebra $A$ corresponds to a~Severi--Brauer variety~$V$ over a field $\mathbf{F}$ then we have~\mbox{$\mathrm{Aut}(V)  \simeq A^*/\mathbf{F}^*.$}
\end{Th}

Theorem~\ref{csaprop}  allows us to obtain restrictions on the orders of automorphisms of Severi--Brauer varieties.

\begin{Lemma}[{cf.~\cite[Lemma 5.2]{Shramov1}}]\label{lemma:l<p}
Let~$V$ be a non-trivial Severi--Brauer variety of odd prime degree~$p$ over~$\mathbb{Q}.$ Then $\mathrm{Aut}(V)$ does not contain elements of prime order~\mbox{$l > p.$} 
\end{Lemma}

\begin{proof}
Suppose that there is an element $x \in \mathrm{Aut}(V)$ of prime order~$l > p.$   Let $A$ be the central simple algebra of dimension $p^2$ over $\mathbb{Q}$ corresponding to the variety $V$. By Theorem~\ref{csaprop} there exists an element $\bar{x} \in A^*$ such that $\bar{x} \notin \mathbb{Q}^*$ and $\bar{x}^l=a \in \mathbb{Q}^*.$

Since $V$ is a non-trivial Severi–Brauer variety, the central simple algebra $A$ is a division algebra by Theorem~\ref{CSAcorps}. Thus $x$ generates a subfield $B \subset A$. Since
$$
\dim_{\mathbb{Q}}A=\dim_{B}A \cdot \dim_{\mathbb{Q}}B,
$$
\noindent we obtain $\dim_{\mathbb{Q}}B=p,$ because $\dim_{\mathbb{Q}}B=p^2$ is impossible as $B$ is a commutative algebra while $A$ is not. 

Let $f(t)$ be a minimal polynomial of $\bar{x}$ over $\mathbb{Q}.$ Its degree is equal to $p.$ Moreover, the polynomial $f(t)$ divides $t^l-a.$ In particular, the polynomial $t^l-a$ is reducible over $\mathbb{Q}.$ So by~\cite[Theorem \rom{6}.9.1]{Lang} we have~\mbox{$a=c^l$} for some $c \in \mathbb{Q}^*.$ We get that the polynomial $f(ct)$ divides the $l$-th cyclotomic polynomial $\Phi_l(t)$, but the latter is irreducible over~$\mathbb{Q}$. The assertion of the lemma follows from this contradiction.

\end{proof}

\begin{Lemma}\label{lemma:z/pz^2}
Let $V$ be a non-trivial Severi--Brauer variety of odd prime degree~$p$ over~$\mathbb{Q}.$ Let $G$ be a finite subgroup in $\mathrm{Aut}(V).$ Then $G \subseteq (\mathbb{Z}/p\mathbb{Z})^2.$
\end{Lemma}

\begin{proof}
By~\cite[Theorem 1.4]{Savelyeva} (see also~\cite[Section~1]{Sonina}) we get that 
$$
G \subseteq (\mathbb{Z}/n\mathbb{Z} \rtimes \mathbb{Z}/p\mathbb{Z}) \times \mathbb{Z}/p\mathbb{Z},
$$
where $n$ is a positive integer such that every prime~$q$ that divides $n$ is congruent to~$1$ modulo $p.$ So by Lemma~\ref{lemma:l<p} we get that $n=1.$ Therefore, we get~\mbox{$G \subseteq (\mathbb{Z}/p\mathbb{Z})^2.$}

\end{proof}

Recall the definition of a cyclic algebra which is one of the simplest examples of a central simple algebra.

\begin{definition}
Let $\mathbf{F}$ be a field. Let $\mathbf{F} \subset \mathbf{L}$ be a Galois extension with Galois group isomorphic to~$\mathbb{Z}/n\mathbb{Z}$ for an integer  $n \geqslant 2.$ Consider an element $a \in \mathbf{F}^*$ and a generator $\sigma$ of  $\mathrm{Gal}(\mathbf{L}/\mathbf{F}).$ Then one can associate with $a$ and $\sigma$ a central simple algebra $A$ over $\mathbf{F}$ which is generated by $\mathbf{L}$ and an element $\alpha$ such that
$$
\alpha^n=a \;\; \text{and} \;\; \lambda\alpha=\alpha\sigma(\lambda) \;\; \text{for} \;\; \lambda \in \mathbf{L}.
$$
Such a central simple algebra $A$ is called a \textit{cyclic algebra}.
\end{definition}

The following ancient result about central simple division algebras over~$\mathbb{Q}$ was proven by titans of algebra: Richard Brauer, Emmy Noether and Helmut Hasse in~1932.

\begin{Th}[{\cite[Hauptsatz]{BNH}}]\label{Emmy}
Let $A$ be a central simple algebra over~$\mathbb{Q}.$ Assume that~$A$ is a division algebra. Then $A$ is a cyclic algebra.
\end{Th}

Now we are going to study $p$-subgroups in the automorphism groups of Severi--Brauer varieties of odd prime degree~$p$. First of all, we show that the group $\mathbb{Z}/p\mathbb{Z}$ acts faithfully on every non-trivial Severi--Brauer variety of prime degree~$p$ over~$\mathbb{Q}.$

\begin{Lemma}\label{lemma:autz/pz}
Let $V$ be a non-trivial Severi--Brauer variety of prime degree~$p$ over the field $\mathbb{Q}.$ Then $\mathrm{Aut}(V)$  contains a subgroup isomorphic to~$\mathbb{Z}/p\mathbb{Z}.$ 
\end{Lemma}

\begin{proof}

 
 The Severi--Brauer variety $V$ corresponds to a division algebra $A$ by Theorem~\ref{CSAcorps}, and this algebra is a cyclic algebra by Theorem~\ref{Emmy}. By~\cite[Proposition~2.5.2]{SimpleAlgebras} the algebra $A$ is generated by a Galois extension $\mathbb{Q} \subset \mathbf{L}$ of degree~$p$ and an element~\mbox{$\alpha \in A$} such that $\alpha \notin \mathbb{Q}^*$ and~\mbox{$\alpha^p \in \mathbb{Q}^*.$} Furthermore, one has $\alpha\lambda=\sigma(\lambda)\alpha$ for all $\lambda \in \mathbf{L},$ where $\sigma$ is a generator of the Galois group $\mathrm{Gal}(\mathbf{L}/\mathbb{Q})$. By Theorem~\ref{csaprop} the element~$\alpha$ gives us an automorphism of order~$p.$  Therefore, we obtain
 $$
 \mathbb{Z}/p\mathbb{Z} \subset \mathrm{Aut}(V).
 $$

\end{proof}

\begin{Lemma}\label{lemma:matrixalgebraifscalar}
Let $A$ be a central simple algebra of dimension~$p^2$ over  $\mathbb{Q}.$  Consider the element $x \in A$ such that $x^p=1$ and~\mbox{$x \notin \mathbb{Q}.$} Then $A$ is a matrix algebra. 

\end{Lemma}

\begin{proof}  
Assume that $A$ is not a matrix algebra. Then by Theorem~\ref{CSAcorps} we get that~$A$ is a division algebra. We have
$$
0=x^{p}-1=(x-1)\sum\limits_{i = 0}^{n - 1} x^{i}
$$
Since $x \notin \mathbb{Q}$ and $A$ is a division algebra, we get that $\sum\limits_{i = 0}^{n - 1} x^{i}=0$. So we have an extension $\mathbb{Q} \subset \mathbf{K}$ of degree $p - 1$ in $A$ generated by $x.$ This means that
$$
p^2=\dim_{\mathbb{Q}}A=\dim_{\mathbb{Q}}\mathbf{K} \cdot \dim_{\mathbf{K}}A=(p - 1)\dim_{\mathbf{K}}A,
$$ 
which is a contradiction.

\end{proof}

\begin{Lemma}\label{lemma:csaomega}
Let $A$ be a central simple algebra of dimension $p^2$ over  $\mathbb{Q}.$ Assume that
$$
\widehat{T}=(\mathbb{Z}/p\mathbb{Z})^2 \subset A^*/\mathbb{Q}^*,
$$
and let   $x$ and~$y$ be  generators of $\widehat{T}.$ If $\bar{x}$ and~$\bar{y}$  are any elements in the preimages of $x$ and $y$  under the natural homomorphism~\mbox{$A^* \to A^*/\mathbb{Q}^*,$} respectively, then 
$$
\bar{x}\bar{y}\bar{x}^{-1}=\zeta_p\bar{y},
$$
\noindent where $\zeta_p$ is a non-trivial $p$-th root of unity.
\end{Lemma}

\begin{proof}
We have $xyx^{-1}y^{-1}=1$ in~$A^*/\mathbb{Q}^*,$ because $x$ and $y$ commute with each other. Therefore, we obtain  
\begin{equation}\label{eq:csaomegaa}
\bar{x}\bar{y}\bar{x}^{-1}=a\bar{y}
\end{equation}
\noindent for some $a \in \mathbb{Q}^*.$ As $\bar{x}^p \in \mathbb{Q}^*,$ we get that $\bar{x}^p$ commutes with $\bar{y},$ so by \eqref{eq:csaomegaa}
$$
\bar{y}=\bar{x}^p\bar{y}\bar{x}^{-p}=a^p\bar{y}.
$$
\noindent Thus, we have $a^p=1.$

By contradiction, assume that $a=1.$  So the elements $\bar{x}$ and $\bar{y}$ commute with each other. Assume that $\bar{x}^p$ and $\bar{y}^p$ belong to $(\mathbb{Q}^*)^p.$ So without loss of generality we can take $\bar{x}$ and $\bar{y}$ in the preimages of $x$ and $y,$ respectively, such that $\bar{x}^p=\bar{y}^p=1$ and $\bar{x}$ and $\bar{y}$ also commute with each other. This means that $(\mathbb{Z}/p\mathbb{Z})^2 \subset A^*.$ By Lemma~\ref{lemma:matrixalgebraifscalar} we get that $A$ is a matrix algebra. So by Lemma~\ref{lemma:representationz/pz} we get that  
$$
(\mathbb{Z}/p\mathbb{Z})^2 \not\subset A^*.
$$
Therefore, in this case $a \neq 1.$

Now assume that  $\bar{x}^p \notin (\mathbb{Q}^*)^p.$ Consider the subalgebra $B$ of $A$ generated by the elements
\begin{equation}\label{eq:elementsofB}
1,\bar{x}^i, \bar{y}^{j}, \bar{x}^i\bar{y}^j, \,\, i,j \in \{1, \cdots , p - 1\} .
\end{equation}

 Let us prove that these elements are linearly independent over $\mathbb{Q}.$ First of all, note that the elements $1$, $\bar{x}, \cdots, \bar{x}^{p - 1}$ are linearly independent over $\mathbb{Q}.$ Indeed, let us consider the algebra $B'$ generated by these $p$ elements. Then from the equation
$$
\dim_{B'}A \cdot \dim_{\mathbb{Q}}B'=\dim_{\mathbb{Q}}A
$$
\noindent we get that $\dim_{\mathbb{Q}}B'$ is either $1,$ or $p.$ However, by assumptions the element $\bar{x}$ does not lie in $\mathbb{Q}.$ As
$$
A \supseteq B \supseteq B',
$$
\noindent we get that $\dim_{\mathbb{Q}}B$ is either $p,$ or $p^2.$ So it is enough to prove that $B' \neq B$ to exclude the case $\dim_{\mathbb{Q}}B=p.$ To prove this it is enough to prove that $\bar{y}$ does not lie in $B'.$ Assume the contrary. Then by Lemma~\ref{lemma:fielddegree3} the element $\bar{y}$ is equal to~$c\bar{x}^i$ for $c \in \mathbb{Q}^*$ and $i \in \{0,1, \ldots, p-1\}.$ That is the contradiction with the fact that $x$ and $y$ generate the group $\left(\mathbb{Z}/p\mathbb{Z}\right)^2.$

 This gives us  $\dim_{\mathbb{Q}}B=p^2=\dim_{\mathbb{Q}}A.$ But~$A$ is a central simple algebra and $B$ is a commutative algebra. This contradiction gives us that $a$ is a non-trivial $p$-th root of unity.

\end{proof}

\begin{Lemma}\label{lemma:z/pz^2notaut}
 Let $V$ be a  non-trivial Severi--Brauer variety of odd prime degree~$p$ over $\mathbb{Q}.$  Then~$(\mathbb{Z}/p\mathbb{Z})^2$ is not a subgroup of $\mathrm{Aut}(V)$.
\end{Lemma}

\begin{proof}

 Let $A$ be a central simple algebra corresponding to the Severi--Brauer variety~$V.$ Then by Theorem \ref{csaprop} we obtain~\mbox{$\mathrm{Aut}(V) \simeq A^*/\mathbb{Q}^*.$}  Assume that
$$
(\mathbb{Z}/p\mathbb{Z})^2 \subset A^*/\mathbb{Q}^*.
$$
\noindent  Let~$x$ and $y$ be generators of $(\mathbb{Z}/p\mathbb{Z})^2.$ Let~$\bar{x}$ and~$\bar{y}$ be any elements in the preimages of $x$ and $y$  under the natural homomorphism~$A^* \to A^*/\mathbb{Q}^*,$ respectively.  By Lemma~\ref{lemma:csaomega} we obtain
$$
\bar{x}\bar{y}\bar{x}^{-1}=\zeta_p\bar{y}
$$
\noindent where $\zeta_p$ is a non-trivial  $p$-th root of unity. However, the field $\mathbb{Q}$ does not contain non-trivial $p$-th root of unity. This contradiction gives us 
$$
(\mathbb{Z}/p\mathbb{Z})^2 \not\subset \mathrm{Aut}(V).
$$

\end{proof}

\section{Proof of Theorem \ref{Z/pZ}}\label{section:proofoftheorem}
In this section we prove Theorem \ref{Z/pZ}.

\begin{proof}[Proof of Theorem \ref{Z/pZ}]
Let $G \subset \mathrm{Aut}(V)$ be a finite subgroup of the automorphism group of $V.$ By Lemma~\ref{lemma:z/pz^2} one has $G \subset (\mathbb{Z}/p\mathbb{Z})^2.$ So by Lemma~\ref{lemma:z/pz^2notaut} the group $G$ is contained in $\mathbb{Z}/p\mathbb{Z}$ and by Lemma~\ref{lemma:autz/pz} the group $\mathbb{Z}/p\mathbb{Z}$ is contained in $\mathrm{Aut}(V)$.

\end{proof}

%
%
%

\providecommand{\bysame}{\leavevmode\hbox to3em{\hrulefill}\thinspace}
\providecommand{\MR}{\relax\ifhmode\unskip\space\fi MR }
\providecommand{\MRhref}[2]{%
  \href{http://www.ams.org/mathscinet-getitem?mr=#1}{#2}
}
\providecommand{\href}[2]{#2}

\end{document}